\documentclass[12pt]{amsart}
\usepackage[english]{babel}
\usepackage{amssymb}

\usepackage{xcolor}
\definecolor{darkgreen}{rgb}{0,0.5,0}
\usepackage[colorlinks=true, linkcolor=darkgreen, 
    urlcolor=blue, citecolor=red]{hyperref}
\usepackage{mathtools}
\mathtoolsset{showonlyrefs}
\usepackage{tikz}
\usetikzlibrary{arrows.meta}

\newtheorem{proposition}{Proposition}
\newtheorem{theorem}{Theorem}
\newtheorem{corollary}{Corollary}
\theoremstyle{definition}
\newtheorem{example}{Example}
\theoremstyle{remark}
\newtheorem{remark}{Remark}

\title{Graph Lie algebras and their central extensions}

\author{Ivan Sechin}
\address{I.S.: BIMSA, Beijing, China}
\email{sechin@bimsa.cn}

\author{Nicolai Reshetikhin}
\address{N.R.: YMSC, Tsinghua University, Beijing; BIMSA, Beijing, China;
Saint Petersburg University, Saint Petersburg, Russia; Department of
Mathematics, University of California, Berkeley}
\email{reshetik@math.berkeley.edu}

\begin{document}

\begin{abstract}
We propose a new class of infinite-dimensional Lie algebras which naturally
generalizes affine Kac--Moody algebras. They are extensions of graph Lie
algebras with special isotropy conditions.
\end{abstract}

\maketitle

\section{Introduction}

Affine Kac--Moody algebras \cite{KacBook} play a fundamental role in modern
mathematical physics, most notably in two-dimensional conformal field
theories and, in particular, in the WZW model \cite{KZ}; see also
\cite{DMS} and the references therein. General Kac--Moody algebras have many
remarkable features, but, with a few exceptions such as the hyperbolic
algebras whose Weyl groups are arithmetic and which are related to
automorphic forms \cite{FF, Borcherds}, they remain largely disconnected
from geometry and physics. Important results were obtained in constructing
infinite-dimensional Lie algebras with a natural geometric origin. Most
notable among them is the algebro-geometric construction of Maulik and
Okounkov \cite{MO}. A differential-geometric construction of central
extensions of current algebras on Riemann surfaces was given by Etingof
and Frenkel \cite{EF}.

In this note we propose a different generalization of affine Kac--Moody
algebras, in which the circle (the ``loop'' in the ``loop algebra'') is
replaced by a finite oriented graph \(\Gamma\). Below we outline the main
results.

Denote the set of vertices of \(\Gamma\) by \(V\) and the set of edges by
\(E\). Assign a finite-dimensional complex simple Lie algebra
\(\mathfrak{g}_e\) to each edge \(e\), and define the Lie algebra
\begin{equation}
    \mathfrak{g}(\Gamma) = \bigoplus_{e \in E} \mathfrak{g}(e).
\end{equation}
Here \(\mathfrak{g}(e)\) is the space of smooth maps
\([0, 1] \to \mathfrak{g}_e\), with derivatives of every order extending
continuously to both endpoints.

For each vertex \(v \in V\), let \(\operatorname{Star}(v)\) be the set of
incident half-edges \((e, i)\), with \(i = 0\) at the source and
\(i = 1\) at the target of \(e\). A loop at \(v\) contributes both
\((e, 0)\) and \((e, 1)\). Define the Lie algebra
\begin{equation}
    \mathfrak{g}_v
        = \bigoplus_{(e, i) \in \operatorname{Star}(v)}
            \mathfrak{g}_e.
\end{equation}

Denote by \(\langle \cdot, \cdot \rangle_e\) the Killing form on
\(\mathfrak{g}_e\). Assign a weight \(k_e \in \mathbb{C}\) to each edge
\(e \in E\) and define the bilinear form
\(\langle \cdot, \cdot \rangle_v\) on \(\mathfrak{g}_v\) as
\begin{equation}\label{eq:form}
    \langle x, y \rangle_v
        = \sum_{(e, 1) \in \operatorname{Star}(v)}
            k_e \langle x_{e, 1}, y_{e, 1} \rangle_e
        - \sum_{(e, 0) \in \operatorname{Star}(v)}
            k_e \langle x_{e, 0}, y_{e, 0} \rangle_e.
\end{equation}
Here \(x, y \in \mathfrak{g}_v\), and \(x_{e, i}, y_{e, i}\) are their
components in the summand indexed by the half-edge \((e, i)\).

For each \(v \in V\), choose a Lie subalgebra
\(\mathfrak{h}_v \subset \mathfrak{g}_v\) that is isotropic with respect
to the form \eqref{eq:form}. Write
\(\mathfrak{h} = (\mathfrak{h}_v)_{v \in V}\). Define the evaluation map
\(\mathrm{ev}_v \colon \mathfrak{g}(\Gamma) \to \mathfrak{g}_v\) by
\begin{equation}
    \gamma_v := \mathrm{ev}_v(\gamma)
        = \big(\gamma_e(i)\big)_{(e, i) \in \operatorname{Star}(v)}.
\end{equation}
Define the Lie subalgebra
\begin{equation}
    \mathfrak{g}(\Gamma, \mathfrak{h})
        = \bigl\{\gamma \in \mathfrak{g}(\Gamma) \;\bigm|\;
            \gamma_v \in \mathfrak{h}_v \subset \mathfrak{g}_v
            \text{ for each } v \in V\bigr\}
        \subset \mathfrak{g}(\Gamma).
\end{equation}

Our main results can be formulated as the following two theorems.

\begin{theorem}
\begin{enumerate}
    \item[(i)] The bilinear form
    \(\omega \colon \mathfrak{g}(\Gamma, \mathfrak{h})
        \times \mathfrak{g}(\Gamma, \mathfrak{h}) \to \mathbb{C}\),
    \begin{equation}
        \omega(\gamma, \eta) = \sum_{e \in E} k_e \int_0^1
            \langle \gamma_e(\tau_e),
                \partial_{\tau_e}\eta_e(\tau_e) \rangle_e
            \mathrm{d}\tau_e,
    \end{equation}
    is a \(2\)-cocycle on \(\mathfrak{g}(\Gamma, \mathfrak{h})\) (in the
    Chevalley--Eilenberg complex with trivial coefficients). This gives the
    central extension \(\widehat{\mathfrak{g}}(\Gamma, \mathfrak{h})\) of
    \(\mathfrak{g}(\Gamma, \mathfrak{h})\):
    \begin{equation}
        \widehat{\mathfrak{g}}(\Gamma, \mathfrak{h})
            = \mathfrak{g}(\Gamma, \mathfrak{h}) \oplus \mathbb{C}K
    \end{equation}
    with the Lie brackets
    \begin{equation}
        [\gamma, \eta]^{\wedge}
            = [\gamma, \eta] + K\,\omega(\gamma, \eta),
        \qquad
        [\gamma, K]^{\wedge} = 0.
    \end{equation}

    \item[(ii)] Fix the parametrizations \(\tau_e \in [0, 1]\) of all
    edges and define \(D\) on \(\mathfrak{g}(\Gamma)\) by
    \begin{equation}
        (D\gamma)_e = \partial_{\tau_e}\gamma_e.
    \end{equation}
    The subalgebra
    \begin{equation}
        \mathfrak{g}_{\infty}(\Gamma, \mathfrak{h})
            = \{\gamma \in \mathfrak{g}(\Gamma) \colon
                \mathrm{ev}_v(D^n\gamma) \in \mathfrak{h}_v
                \text{ for all }v \in V\text{ and }n \geq 0\}
    \end{equation}
    is stable under \(D\). This derivation lifts to the restricted
    central extension
    \begin{equation}
        \widehat{\mathfrak{g}}_{\infty}(\Gamma, \mathfrak{h})
            = \mathfrak{g}_{\infty}(\Gamma, \mathfrak{h})
                \oplus \mathbb{C}K
    \end{equation}
    by setting \(D(K) = 0\). Its extension by this derivation is the
    semidirect product
    \begin{equation}
        \widetilde{\mathfrak{g}}_{\infty}(\Gamma, \mathfrak{h})
            = \widehat{\mathfrak{g}}_{\infty}(\Gamma, \mathfrak{h})
                \rtimes \mathbb{C}d,
    \end{equation}
    where
    \begin{equation}
        [d, \gamma] = D\gamma,
        \qquad
        [d, K] = 0.
    \end{equation}

    \item[(iii)] The Lie algebra
    \(\widetilde{\mathfrak{g}}_{\infty}(\Gamma, \mathfrak{h})\) has a
    symmetric invariant bilinear form \(B\). For
    \(\gamma, \eta \in \mathfrak{g}_{\infty}(\Gamma, \mathfrak{h})\),
    its restriction to currents is given by
    \begin{equation}
        B(\gamma, \eta) = B_0(\gamma, \eta)
            = \sum_{e \in E} k_e \int_0^1
                \langle \gamma_e(\tau_e), \eta_e(\tau_e) \rangle_e
                \mathrm{d}\tau_e,
    \end{equation}
    and the remaining pairings are
    \begin{gather}
        B(\gamma, K) = B(\gamma, d) = 0,
        \qquad
        B(K, d) = -1,
        \qquad
        B(K, K) = B(d, d) = 0.
    \end{gather}
    The form \(B\) is nondegenerate if and only if \(B_0\) is
    nondegenerate on \(\mathfrak{g}_{\infty}(\Gamma, \mathfrak{h})\).
    These equivalent conditions hold if and only if \(k_e \neq 0\) for
    every edge \(e \in E\).
\end{enumerate}
\end{theorem}

Thus, \(\widetilde{\mathfrak{g}}_{\infty}(\Gamma, \mathfrak{h})\) can be
regarded as a generalization of affine Kac--Moody algebras.

Fix the edge parametrizations \(\tau_e\). The \emph{Lie algebra of graph
vector fields} \(\operatorname{Vec}(\Gamma)\) consists of the edgewise
vector fields
\begin{equation}
    X_f = \big(f_e(\tau_e)\partial_{\tau_e}\big)_{e \in E},
    \qquad f_e \in C^\infty([0, 1], \mathbb{C}),
\end{equation}
satisfying the vertex conditions
\begin{equation}
    \partial_{\tau_e}^{r}f_e(i)
        = \partial_{\tau_{e'}}^{r}f_{e'}(j)
    \quad\text{for }(e, i), (e', j) \in \operatorname{Star}(v),
    \quad v \in V,\quad r \in \mathbb{Z}_{\geq 0}.
\end{equation}
Here \(i = 1\) for an incoming half-edge and \(i = 0\) for an outgoing
half-edge; both half-edges of a loop are included. The Lie bracket is
given edgewise by
\begin{equation}
    [X_f, X_g]_e
        = \big(f_e\partial_{\tau_e}g_e
            - g_e\partial_{\tau_e}f_e\big)\partial_{\tau_e}.
\end{equation}
Define the graph analogue of the Gelfand--Fuchs bilinear form by
\begin{equation}
    \omega_{\mathrm{GF}}(X_f, X_g)
        = \sum_{e \in E} \int_0^1
            \big(g_e\partial_{\tau_e}^{3}f_e
                - f_e\partial_{\tau_e}^{3}g_e\big)\,\mathrm{d}\tau_e,
    \qquad X_f, X_g \in \operatorname{Vec}(\Gamma).
\end{equation}

\begin{theorem}\label{thm:graph-virasoro}
\begin{enumerate}
    \item[(i)] If every vertex of \(\Gamma\) has equally many incoming
    and outgoing half-edges, then \(\omega_{\mathrm{GF}}\) is a
    two-cocycle on \(\operatorname{Vec}(\Gamma)\) with trivial
    coefficients.

    \item[(ii)] Under this assumption, the corresponding central
    extension, the \emph{graph Virasoro algebra}
    \begin{equation}
        \operatorname{Vir}(\Gamma)
            = \operatorname{Vec}(\Gamma) \oplus \mathbb{C}c,
    \end{equation}
    with brackets
    \begin{equation}
        [X_f, X_g]_{\operatorname{Vir}}
            = [X_f, X_g] + \omega_{\mathrm{GF}}(X_f, X_g)c,
        \qquad [X_f, c]_{\operatorname{Vir}} = 0,
    \end{equation}
    acts on
    \(\widehat{\mathfrak{g}}_{\infty}(\Gamma, \mathfrak{h})\)
    by derivations. This defines the semidirect-product extension
    \begin{equation}
        \widetilde{\mathfrak{g}}^{\operatorname{Vir}}_{\infty}
            (\Gamma, \mathfrak{h})
            = \widehat{\mathfrak{g}}_{\infty}(\Gamma, \mathfrak{h})
                \rtimes \operatorname{Vir}(\Gamma),
    \end{equation}
    whose mixed brackets, for
    \(\gamma \in \mathfrak{g}_{\infty}(\Gamma, \mathfrak{h})\), are
    \begin{gather}
        [X_f, \gamma] = D_{X_f}\gamma,
        \qquad (D_{X_f}\gamma)_e
            = f_e\partial_{\tau_e}\gamma_e, \\
        [X_f, K] = 0,
        \qquad [c, \gamma] = 0,
        \qquad [c, K] = 0.
    \end{gather}
\end{enumerate}
\end{theorem}

The structure of the paper is as follows. In Section~2 we introduce graph
Lie algebras and describe their basic properties. Section~3 contains the
proof of the cocycle property of \(\omega\) and constructs the extension
\(\widetilde{\mathfrak{g}}_{\infty}(\Gamma, \mathfrak{h})\) by a
derivation in Subsection~\ref{subsec:extensions}.
Subsection~\ref{subsec:invariant-form} proves the invariance of its
bilinear form and characterizes its nondegeneracy.
Section~\ref{sec:vector-fields} introduces graph vector fields and the
graph Virasoro algebra, together with their actions on the centrally
extended current algebra, and proves Theorem~\ref{thm:graph-virasoro}.

\subsection{Acknowledgements}

The authors are grateful to I.~Frenkel for inspiring discussions. The
research of N.R. was supported by the Changjiang fund, by the Collaboration
Grant ``Categorical Symmetries'' from the Simons Foundation, and by the
project 075-15-2025-013 funded by the Ministry of Science and Higher
Education of the Russian Federation.

\section{Graph Lie algebras}

\subsection{Graphs and edge Lie algebras}

Let \(\Gamma\) be a finite oriented graph with possible loops and multiple
edges. It consists of a finite set of vertices \(V\), a finite set of edges
\(E\), and an incidence map
\begin{equation}
    (s, t) \colon E \to V \times V
    \qquad
    \begin{tikzpicture}[baseline=-0.5ex]
        \fill (0, 0) circle (1.5pt)
            node[below] {\(s(e)\)};
        \fill (2.5, 0) circle (1.5pt)
            node[below] {\(t(e)\)};
        \draw[->] (0.1, 0) -- (2.4, 0)
            node[midway, above] {\(e\)};
    \end{tikzpicture}
\end{equation}
For an edge \(e \in E\), the vertices \(s(e)\) and \(t(e)\) are called its 
source and target, respectively. Edges with \(s(e) = t(e)\) are called loops. 
Since the map \((s, t)\) is not required to be injective, multiple edges 
between two vertices are allowed.

We identify each edge \(e \in E\) with an interval in a way compatible with
the orientation. We choose local parametrizations of the edges
\(\tau_e \in [0, 1]\) such that \(\tau_e = 0\) for the source vertex
\(s(e)\) and \(\tau_e = 1\) for the target vertex \(t(e)\).
Let \(\{\mathfrak{g}_e\}_{e \in E}\) be a collection of finite-dimensional
complex simple Lie algebras indexed by edges of \(\Gamma\).

For each edge \(e\) define a Lie algebra \(\mathfrak{g}(e)\) as the space of
smooth maps
\begin{equation}
    \gamma_e \colon [0, 1] \to \mathfrak{g}_e.
\end{equation}
Here smoothness means infinite differentiability on \((0, 1)\), with the
map and every derivative extending continuously to \([0, 1]\). Derivatives
at \(0\) and \(1\) are the corresponding one-sided limits. No analyticity
is assumed.
It has the pointwise Lie bracket
\begin{equation}
    [\gamma_e, \eta_e](\tau_e)
        = [\gamma_e(\tau_e), \eta_e(\tau_e)]_{\mathfrak{g}_e}.
\end{equation}

Define the graph Lie algebra \(\mathfrak{g}(\Gamma)\) as a direct sum
\begin{equation}
    \mathfrak{g}(\Gamma) = \bigoplus_{e \in E} \mathfrak{g}(e),
\end{equation}
with the pointwise Lie bracket on each summand.

\subsection{Vertex conditions}

For each vertex \(v \in V\), let
\begin{equation}
    \operatorname{Star}(v)
        = \{(e, 0) \in E \times \{0\} \colon s(e) = v\}
            \sqcup
          \{(e, 1) \in E \times \{1\} \colon t(e) = v\}.
\end{equation}
This is the set of half-edges incident to \(v\). Define
\begin{equation}
    \mathfrak{g}_v
        = \bigoplus_{(e, i) \in \operatorname{Star}(v)}
            \mathfrak{g}_e.
\end{equation}
The evaluation map
\(\mathrm{ev}_v \colon \mathfrak{g}(\Gamma) \to \mathfrak{g}_v\) is
\begin{equation}
    \gamma_v := \mathrm{ev}_v(\gamma)
        = \big(\gamma_e(i)\big)_{(e, i) \in \operatorname{Star}(v)}.
\end{equation}
For a loop based at \(v\), the two half-edges \((e, 0)\) and \((e, 1)\)
are distinct, so both endpoint values occur in \(\gamma_v\).

Denote by \(\langle \cdot, \cdot \rangle_e\) the Killing form on 
\(\mathfrak{g}_e\). Fix the collection of complex numbers \(\{k_e\}_{e \in E}\)
indexed by edges of \(\Gamma\). For each vertex \(v \in V\), define
a weighted bilinear form on \(\mathfrak{g}_v\) by
\begin{equation}
    \langle \gamma_v, \eta_v \rangle_v = 
        \sum_{e \colon v = t(e)} 
            k_e \langle \gamma_e(1), \eta_e(1) \rangle_e -
        \sum_{e \colon v = s(e)}
            k_e \langle \gamma_e(0), \eta_e(0) \rangle_e.
\end{equation}
Since each \(\mathfrak{g}_e\) is simple and the Killing form is invariant on
\(\mathfrak{g}_e\), \(\langle \cdot, \cdot \rangle_v\) is an invariant
bilinear form on \(\mathfrak{g}_v\).

For each vertex \(v \in V\), fix a Lie subalgebra
\(\mathfrak{h}_v \subset \mathfrak{g}_v\) isotropic with respect to
\(\langle \cdot, \cdot \rangle_v\), i.e.,
\begin{equation}
    \langle \gamma_v, \eta_v \rangle_v = 0, \quad
        \text{for all } \gamma_v, \eta_v \in \mathfrak{h}_v. 
\end{equation}
Write \(\mathfrak{h} = (\mathfrak{h}_v)_{v \in V}\). Define
\(\mathfrak{g}(\Gamma, \mathfrak{h}) \subset \mathfrak{g}(\Gamma)\) to be the
Lie subalgebra consisting of the elements
\(\gamma = (\gamma_e)_{e \in E} \in \mathfrak{g}(\Gamma)\) whose vertex
values satisfy \(\gamma_v \in \mathfrak{h}_v\) for all \(v \in V\).

\subsection{Connected graph Lie algebras}\label{subsec:connected}

The algebra \(\mathfrak{g}(\Gamma)\) is the direct sum of the edge algebras
\(\mathfrak{g}(e)\). The vertex conditions defining
\(\mathfrak{g}(\Gamma, \mathfrak{h})\), however, may couple summands
corresponding to different edges. We describe this coupling by a cutting
operation.

A vertex \(v\) is called \emph{splittable} if there is a partition of its
incident half-edges
\begin{equation}
    \operatorname{Star}(v) = A \sqcup B,
    \qquad
    A \neq \varnothing,
    \qquad
    B \neq \varnothing,
\end{equation}
and Lie subalgebras \(\mathfrak{h}_{v, A}\) and
\(\mathfrak{h}_{v, B}\) such that the isotropic subalgebra
\(\mathfrak{h}_v\) decomposes as a direct sum
\begin{equation}
    \mathfrak{h}_v
        = \mathfrak{h}_{v, A} \oplus \mathfrak{h}_{v, B},
\end{equation}
Here \(\mathfrak{h}_{v, A}\) lies in the direct sum of the edge Lie
algebras corresponding to the half-edges in \(A\), while
\(\mathfrak{h}_{v, B}\) lies in the analogous direct sum for \(B\):
\begin{equation}
    \mathfrak{h}_{v, A}
        \subseteq \bigoplus_{(e, i) \in A} \mathfrak{g}_e,
    \qquad
    \mathfrak{h}_{v, B}
        \subseteq \bigoplus_{(e, i) \in B} \mathfrak{g}_e.
\end{equation}
Such a partition is called a \emph{splitting} of \(v\).

Suppose that \(v\) is splittable. Define a new graph \(\Gamma'\) by replacing
\(v\) with two vertices \(v_A\) and \(v_B\). Attach a half-edge to \(v_A\)
if it belongs to \(A\), and to \(v_B\) if it belongs to \(B\). All other
vertices, edges, and incidence relations remain unchanged. Set
\begin{align}
    \mathfrak{h}'_{v_A}
        &= \mathfrak{h}_{v, A},
        \\
    \mathfrak{h}'_{v_B}
        &= \mathfrak{h}_{v, B},
\end{align}
and let \(\mathfrak{h}'_u = \mathfrak{h}_u\) for every other vertex \(u\).
The new vertex subalgebras are isotropic for the restrictions of the
corresponding vertex forms.

Cutting a splittable vertex does not change the associated graph Lie algebra.
More precisely, the identity on the edge algebras induces a canonical
isomorphism
\begin{equation}
    \mathfrak{g}(\Gamma, \mathfrak{h})
        \cong \mathfrak{g}(\Gamma', \mathfrak{h}').
\end{equation}

The graph Lie algebra \(\mathfrak{g}(\Gamma, \mathfrak{h})\) is
\emph{connected} if \(\Gamma\) is connected and no vertex is splittable.

Every graph Lie algebra associated with a finite graph admits, after finitely
many vertex splittings, a decomposition
\begin{equation}\label{eq:connected-decomposition}
    \mathfrak{g}(\Gamma, \mathfrak{h})
        \cong \bigoplus_{j = 1}^r
            \mathfrak{g}(\Gamma_j, \mathfrak{h}^{(j)}),
\end{equation}
where every graph Lie algebra
\(\mathfrak{g}(\Gamma_j, \mathfrak{h}^{(j)})\) is connected.

The isotropy conditions at the vertices couple the Lie algebras associated
with individual edges. After the splitting procedure, they therefore give a
direct-sum decomposition into connected graph Lie algebras, rather than into
the Lie algebras associated with individual edges.

\subsection{Subalgebras}

Let \(\widetilde{\Gamma} = (\widetilde V, \widetilde E)\) be a subgraph of
\(\Gamma\). For every \(v \in \widetilde V\), define the
\(\widetilde{\Gamma}\)-supported summand of \(\mathfrak{g}_v\) by
\begin{equation}
    \mathfrak{g}_v^{\widetilde{\Gamma}}
        = \bigoplus_{\substack{
            (e, i) \in \operatorname{Star}(v) \\
            e \in \widetilde E
        }} \mathfrak{g}_e
        \subseteq \mathfrak{g}_v.
\end{equation}
Choose a Lie subalgebra
\begin{equation}
    \widetilde{\mathfrak{h}}_v
        \subseteq \mathfrak{h}_v
            \cap \mathfrak{g}_v^{\widetilde{\Gamma}}
\end{equation}
at every \(v \in \widetilde V\), and write
\(\widetilde{\mathfrak{h}}
    = (\widetilde{\mathfrak{h}}_v)_{v \in \widetilde V}\).
Each \(\widetilde{\mathfrak{h}}_v\) is isotropic for the restriction of the
vertex form to \(\mathfrak{g}_v^{\widetilde{\Gamma}}\). Thus the edge Lie
algebras and weights restricted to \(\widetilde{\Gamma}\), together with
\(\widetilde{\mathfrak{h}}\), define the graph Lie algebra
\(\mathfrak{g}(\widetilde{\Gamma}, \widetilde{\mathfrak{h}})\).

After extending its elements by zero on the edges in
\(E \setminus \widetilde E\), we identify this graph Lie algebra with the
subalgebra
\begin{equation}
    \mathfrak{g}(\widetilde{\Gamma}, \widetilde{\mathfrak{h}})
        = \left\{\gamma \in \mathfrak{g}(\Gamma, \mathfrak{h})
        \;\middle|\;
        \begin{aligned}
            &\gamma_e = 0
                &&\text{for every } e \notin \widetilde E, \\
            &\gamma_v \in \widetilde{\mathfrak{h}}_v
                &&\text{for every } v \in \widetilde V
        \end{aligned}
        \right\}
        \subseteq \mathfrak{g}(\Gamma, \mathfrak{h}).
\end{equation}

\subsection{Basic examples}

\begin{example}
Let \(\Gamma\) have vertices \(v_0, v_1\) and one edge \(e\),
oriented from \(v_0\) to \(v_1\):
\begin{center}
    \begin{tikzpicture}
        \fill (0, 0) circle (1.5pt)
            node[below] {\(v_0\)};
        \fill (2.5, 0) circle (1.5pt)
            node[below] {\(v_1\)};
        \draw[->] (0.1, 0) -- (2.4, 0)
            node[midway, above] {\(e\)};
    \end{tikzpicture}
\end{center}
Set \(\mathfrak{g}_e = \mathfrak{g}\) and \(k_e = k\), where
\(\mathfrak{g}\) is simple. Then
\begin{equation}
    \mathfrak{g}(\Gamma) = \mathfrak{g}(e).
\end{equation}
The vertex Lie algebras are
\(\mathfrak{g}_{v_0} = \mathfrak{g}_{v_1} = \mathfrak{g}\), with forms
\begin{equation}
    \langle x, y \rangle_{v_0} = -k \langle x, y \rangle_e,
    \qquad
    \langle x, y \rangle_{v_1} = k \langle x, y \rangle_e.
\end{equation}
For isotropic Lie subalgebras
\(\mathfrak{h}_{v_0}, \mathfrak{h}_{v_1} \subset \mathfrak{g}\), we obtain
\begin{equation}
    \mathfrak{g}(\Gamma, \mathfrak{h})
        = \{\gamma \in \mathfrak{g}(e) \colon
            \gamma(0) \in \mathfrak{h}_{v_0},
            \ \gamma(1) \in \mathfrak{h}_{v_1}\}.
\end{equation}
In particular, taking
\(\mathfrak{h}_{v_0} = \mathfrak{h}_{v_1} = 0\) gives the algebra of
smooth paths vanishing at both endpoints.
\end{example}

For example, after choosing Borel subalgebras of \(\mathfrak{g}\), one can
take \(\mathfrak{h}_{v_0}\) and \(\mathfrak{h}_{v_1}\) to be their nilpotent
radicals. These subalgebras are isotropic for the Killing form.

\medskip

\begin{example}
Let \(\Gamma\) have one vertex \(v\) and one loop \(e\):
\begin{center}
    \begin{tikzpicture}
        \draw[->] (1, 0) arc (0:350:1)
            node[pos=0.5, left] {\(e\)};
        \fill (1, 0) circle (1.5pt)
            node[right] {\(v\)};
    \end{tikzpicture}
\end{center}
Set \(\mathfrak{g}_e = \mathfrak{g}\) and \(k_e = k\), where
\(\mathfrak{g}\) is simple. Here
\(\mathfrak{g}_v = \mathfrak{g} \oplus \mathfrak{g}\), and its form is
\begin{equation}
    \langle (x_0, x_1), (y_0, y_1) \rangle_v
        = k \big(\langle x_1, y_1 \rangle_e
            - \langle x_0, y_0 \rangle_e\big).
\end{equation}
Take the diagonal subalgebra
\begin{equation}
    \mathfrak{h}_v = \Delta \mathfrak{g}
        = \{(x, x) \colon x \in \mathfrak{g}\}
        \subset \mathfrak{g}_v.
\end{equation}
It is isotropic. The corresponding graph Lie algebra is
\begin{equation}
    \mathfrak{g}(\Gamma, \mathfrak{h})
        = \{\gamma \in \mathfrak{g}(e) \colon
            \gamma(0) = \gamma(1)\}.
\end{equation}
\end{example}

\begin{remark}
The condition defining \(\mathfrak{g}(\Gamma, \mathfrak{h})\) identifies only
the endpoint values \(\gamma(0)\) and \(\gamma(1)\). The smooth loop algebra
\(\mathcal{L}\mathfrak{g} = C^\infty(S^1, \mathfrak{g})\) is therefore
a Lie subalgebra of
\(\mathfrak{g}(\Gamma, \mathfrak{h})\). In the chosen parametrization, it is
given by
\begin{equation}
    \mathcal{L}\mathfrak{g}
        = \{\gamma \in \mathfrak{g}(e) \colon
            \partial_\tau^n \gamma(0)
                = \partial_\tau^n \gamma(1)
            \text{ for all } n \geq 0\}
        \subsetneq \mathfrak{g}(\Gamma, \mathfrak{h}).
\end{equation}
The inclusion is proper because the graph Lie algebra imposes no matching
conditions on the higher endpoint derivatives. For example, for
\(0 \neq x \in \mathfrak{g}\), the path
\begin{equation}
    \gamma(\tau) = \sin(\pi \tau) \cdot x
\end{equation}
satisfies \(\gamma(0) = \gamma(1) = 0\), but
\begin{equation}
    \partial_\tau \gamma(0) = \pi x,
    \qquad
    \partial_\tau \gamma(1) = -\pi x.
\end{equation}
In particular,
\(\partial_\tau \gamma(0) \neq \partial_\tau \gamma(1)\).
Thus \(\gamma\) belongs to \(\mathfrak{g}(\Gamma, \mathfrak{h})\), but not to
\(\mathcal{L}\mathfrak{g}\). Hence the graph Lie algebra is substantially
larger than the smooth loop algebra. The latter contains the real-analytic
loop algebra and the Laurent-polynomial loop algebra used in the
algebraic construction of affine Kac--Moody algebras.
\end{remark}

\section{Extensions of graph algebras}

Following the decomposition in Subsection~\ref{subsec:connected}, we first
assume that \(\mathfrak{g}(\Gamma, \mathfrak{h})\) is connected. We return
to an arbitrary graph Lie algebra after constructing the canonical
edge-integral two-cocycle.

\subsection{The canonical edge-integral two-cocycle}

Define a bilinear map
\begin{equation}
    \omega \colon
        \mathfrak{g}(\Gamma, \mathfrak{h}) \times
        \mathfrak{g}(\Gamma, \mathfrak{h})
        \longrightarrow \mathbb{C}
\end{equation}
by
\begin{equation}
    \omega(\gamma, \eta)
        = \sum_{e \in E} k_e \int_0^1
            \langle \gamma_e, \mathrm{d}\eta_e \rangle_e
        = \sum_{e \in E} k_e \int_0^1
            \langle \gamma_e,
                \partial_{\tau_e}\eta_e \rangle_e\,\mathrm{d}\tau_e.
\end{equation}
Here the arguments \(\tau_e\) of the paths are suppressed.

The edge-integral expression is independent of orientation-preserving
parametrizations. Indeed, let \(\tau_e = \phi_e(\sigma_e)\), where
\(\phi_e \colon [0, 1] \to [0, 1]\) is an orientation-preserving smooth
diffeomorphism, so that \(\phi_e(0) = 0\) and \(\phi_e(1) = 1\). Write
\(\widetilde{\gamma}_e = \gamma_e \circ \phi_e\) and
\(\widetilde{\eta}_e = \eta_e \circ \phi_e\). The chain rule and change
of variables give
\begin{align}
    \int_0^1 \langle \widetilde{\gamma}_e,
        \partial_{\sigma_e}\widetilde{\eta}_e \rangle_e
        \,\mathrm{d}\sigma_e
        &= \int_0^1 \langle \gamma_e \circ \phi_e,
            (\partial_{\tau_e}\eta_e) \circ \phi_e \rangle_e
            \phi'_e\,\mathrm{d}\sigma_e \notag\\
        &= \int_0^1 \langle \gamma_e,
            \partial_{\tau_e}\eta_e \rangle_e\,\mathrm{d}\tau_e.
\end{align}
Thus each edge integral, and hence \(\omega\), is well defined
independently of these parametrizations. Smooth reparametrizations also
preserve the path space and the endpoint values, hence the vertex
conditions.

\begin{proposition}\label{prop:two-cocycle}
The bilinear map \(\omega\) is a Chevalley--Eilenberg two-cocycle on
\(\mathfrak{g}(\Gamma, \mathfrak{h})\) with coefficients in the trivial
\(\mathfrak{g}(\Gamma, \mathfrak{h})\)-module \(\mathbb{C}\).
\end{proposition}

\begin{proof}
First, integration by parts gives
\begin{align}
    \omega(\gamma, \eta) + \omega(\eta, \gamma)
        &= \sum_{e \in E} k_e \int_0^1 \mathrm{d}
            \langle \gamma_e, \eta_e \rangle_e \notag\\
        &= \sum_{e \in E} k_e
            \big(\langle \gamma_e(1), \eta_e(1) \rangle_e
                - \langle \gamma_e(0), \eta_e(0) \rangle_e\big)
                \notag\\
        &= \sum_{v \in V}
            \langle \gamma_v, \eta_v \rangle_v = 0.
\end{align}
The last equality follows from
\(\gamma_v, \eta_v \in \mathfrak{h}_v\) and the isotropy of
\(\mathfrak{h}_v\). Hence \(\omega\) is skew-symmetric.

It remains to verify the cocycle condition. By invariance of the Killing
form, on every edge one has
\begin{equation}
    \langle [\gamma_e, \eta_e], \mathrm{d}\lambda_e \rangle_e
        + \langle [\eta_e, \lambda_e], \mathrm{d}\gamma_e \rangle_e
        + \langle [\lambda_e, \gamma_e], \mathrm{d}\eta_e \rangle_e
        = \mathrm{d}
        \langle [\gamma_e, \eta_e], \lambda_e \rangle_e.
\end{equation}
Therefore
\begin{equation}
    \omega([\gamma, \eta], \lambda)
        + \omega([\eta, \lambda], \gamma)
        + \omega([\lambda, \gamma], \eta)
        = \sum_{v \in V}
        \langle [\gamma_v, \eta_v], \lambda_v \rangle_v = 0.
\end{equation}
Indeed, \([\gamma_v, \eta_v] \in \mathfrak{h}_v\), because
\(\mathfrak{h}_v\) is a Lie subalgebra, and \(\mathfrak{h}_v\) is isotropic.
\end{proof}

\begin{remark}
Assume that \(k_e \neq 0\) for at least one edge \(e\). Then \(\omega\) is
not a coboundary. Choose \(x \in \mathfrak{g}_e\) with
\(\langle x, x \rangle_e \neq 0\), a nonzero real-valued function
\(\rho \in C_c^\infty((0, 1), \mathbb{R})\), and set
\begin{equation}
    f(\tau) = \rho(\tau),
    \qquad
    g(\tau) = \tau\rho(\tau).
\end{equation}
Let \(\gamma\) and \(\eta\) vanish on every other edge, and set
\begin{equation}
    \gamma_e(\tau_e) = f(\tau_e)x,
    \qquad
    \eta_e(\tau_e) = g(\tau_e)x.
\end{equation}
All their endpoint derivatives vanish, so they belong to
\(\mathfrak{g}(\Gamma, \mathfrak{h})\), and
\([\gamma, \eta] = 0\). Every two-coboundary vanishes on this pair, whereas
\begin{equation}
    \omega(\gamma, \eta)
        = \frac{k_e}{2}\langle x, x \rangle_e
            \int_0^1 \rho(\tau)^2\,\mathrm{d}\tau \neq 0.
\end{equation}
Thus \([\omega] \neq 0\) in Chevalley--Eilenberg cohomology.
\end{remark}

If the graph Lie algebra is not connected, apply the same construction to
each connected summand in \eqref{eq:connected-decomposition}. Extend the
canonical cocycle of each summand by zero on all the other summands. If
every connected summand contains an edge \(e\) with \(k_e \neq 0\), a
decomposition with \(r\) connected components yields \(r\) nontrivial
canonical cocycles whose cohomology classes are linearly independent.

\subsection{Additional two-cocycles}

Every graph Lie algebra has the canonical edge-integral cocycle
constructed above. The same construction applies to every choice of the
edge Lie algebras, the isotropic subalgebras, and the graph structure,
without requiring special properties of these data. For special choices,
however, the graph Lie algebra can have additional nontrivial two-cocycles.

The canonical cocycle uses the edge weights \(k_e\) appearing in the vertex
forms. We may replace these coefficients by another collection
\(a = (a_e)_{e \in E} \in \mathbb{C}^E\). For
\(x, y \in \mathfrak{h}_v\), write
\begin{equation}
    x = (x_{e, i})_{(e, i) \in \operatorname{Star}(v)},
    \qquad
    y = (y_{e, i})_{(e, i) \in \operatorname{Star}(v)}.
\end{equation}
The same integration-by-parts argument defines a two-cocycle provided that
the boundary terms still vanish, i.e., if
\begin{equation}
    \sum_{(e, 1) \in \operatorname{Star}(v)}
        a_e \langle x_{e, 1}, y_{e, 1} \rangle_e -
    \sum_{(e, 0) \in \operatorname{Star}(v)}
        a_e \langle x_{e, 0}, y_{e, 0} \rangle_e = 0 \qquad
            \text{for all }v \in V\text{ and }x, y \in \mathfrak{h}_v.
\end{equation}
Such a collection is called an \emph{admissible weight}.
Denote the space of admissible weights by \(\mathcal{W}\).
The defining conditions are linear in \(a\), so \(\mathcal{W}\) is a
vector subspace of \(\mathbb{C}^E\).
The original weight \(k = (k_e)_{e \in E}\) belongs to \(\mathcal{W}\).

For every \(a \in \mathcal{W}\), define
\begin{equation}
    \omega_a(\gamma, \eta)
        = \sum_{e \in E} a_e \int_0^1
            \langle \gamma_e(\tau_e),
                \partial_{\tau_e}\eta_e(\tau_e) \rangle_e
            \mathrm{d}\tau_e.
\end{equation}
In particular, \(\omega_k = \omega\).
The proof of Proposition~\ref{prop:two-cocycle} applies without change.
If \(a \neq 0\), the same argument proves that \(\omega_a\) is not a
coboundary.

The space \(\mathcal{W}\) can have dimension greater than one. For
example, let \(\Gamma\) have one vertex and two loop edges, with
\(\mathfrak{g}_{e_1} = \mathfrak{g}_{e_2} = \mathfrak{g}\), and take
\begin{equation}
    \mathfrak{h}_v = \Delta_4\mathfrak{g}
        = \{(x, x, x, x) \colon x \in \mathfrak{g}\}
        \subset \mathfrak{g}^{\oplus 4}.
\end{equation}
The diagonal intersects every proper coordinate summand trivially, so the
vertex is not splittable. Thus this graph Lie algebra is connected, but
the source and target contributions of each loop cancel separately on
\(\Delta_4\mathfrak{g}\). Therefore every
\(a = (a_{e_1}, a_{e_2}) \in \mathbb{C}^2\) is admissible, and
\(\mathcal{W} = \mathbb{C}^2\).

The isotropic Lie subalgebras attached to the vertices can have nontrivial
Chevalley--Eilenberg cohomology of their own. Fix \(v \in V\), and let
\(\beta_v\) be a Chevalley--Eilenberg two-cocycle on
\(\mathfrak{h}_v\) with coefficients in the trivial module \(\mathbb{C}\).
The restriction of \(\mathrm{ev}_v\) to
\(\mathfrak{g}(\Gamma, \mathfrak{h})\) is a Lie algebra homomorphism to
\(\mathfrak{h}_v\). Its pullback is
\begin{equation}
    \omega_{v, \beta_v}(\gamma, \eta)
        = \beta_v\big(
            \mathrm{ev}_v(\gamma), \mathrm{ev}_v(\eta)\big)
\end{equation}
and is therefore a two-cocycle on
\(\mathfrak{g}(\Gamma, \mathfrak{h})\).

For a nontrivial example, let \(\Gamma\) have one edge from \(v_0\) to
\(v_1\), let \(\mathfrak{g}_e = \mathfrak{sl}_3(\mathbb{C})\), set
\(k_e = 1\), and take
\begin{equation}
    \mathfrak{h}_{v_0}
        = \operatorname{span}\{E_{12}, E_{13}\},
    \qquad
    \mathfrak{h}_{v_1} = 0.
\end{equation}
The isotropic subalgebra \(\mathfrak{h}_{v_0}\) is abelian. Hence every
alternating form on it is a Chevalley--Eilenberg two-cocycle. Let
\(\beta_{v_0}\) be determined by
\begin{equation}
    \beta_{v_0}(E_{12}, E_{13}) = 1.
\end{equation}
For the commuting paths
\begin{equation}
    \gamma(\tau) = (1 - \tau)E_{12},
    \qquad
    \eta(\tau) = (1 - \tau)E_{13},
\end{equation}
one has
\begin{equation}
    \omega_{v_0, \beta_{v_0}}(\gamma, \eta) = 1,
    \qquad
    \omega(\gamma, \eta) = 0.
\end{equation}
Every coboundary vanishes on this commuting pair. Hence this cocycle is
nontrivial and its cohomology class is independent of \([\omega]\).

\subsection{Central extension and extension by derivations}
\label{subsec:extensions}

Consider the vector space
\begin{equation}\label{eq:central-extension}
    \widehat{\mathfrak{g}}(\Gamma, \mathfrak{h})
        = \mathfrak{g}(\Gamma, \mathfrak{h}) \oplus \mathbb{C}K.
\end{equation}
For \(a, b \in \mathbb{C}\) and
\(\gamma, \eta \in \mathfrak{g}(\Gamma, \mathfrak{h})\), set
\begin{equation}\label{eq:central-extension-bracket}
    [\gamma + aK, \eta + bK]^{\wedge}
        = [\gamma, \eta] + \omega(\gamma, \eta)K.
\end{equation}

\begin{corollary}
The vector space \(\widehat{\mathfrak{g}}(\Gamma, \mathfrak{h})\), with the
bracket in \eqref{eq:central-extension-bracket}, is a one-dimensional
central extension of \(\mathfrak{g}(\Gamma, \mathfrak{h})\).
\end{corollary}

\begin{proof}
This follows immediately from Proposition~\ref{prop:two-cocycle}.
\end{proof}

Fix the parametrizations \(\tau_e \in [0, 1]\) of all edges \(e \in E\).
Define edgewise differentiation on \(\mathfrak{g}(\Gamma)\) by
\begin{equation}
    (D\gamma)_e = \partial_{\tau_e}\gamma_e.
\end{equation}
Although \(D\) is a derivation of the full graph Lie algebra
\(\mathfrak{g}(\Gamma)\), it need not preserve the vertex conditions defining
\(\mathfrak{g}(\Gamma, \mathfrak{h})\). Consider instead the subalgebra
\begin{equation}
    \mathfrak{g}_{\infty}(\Gamma, \mathfrak{h})
        = \{\gamma \in \mathfrak{g}(\Gamma) \colon
            \mathrm{ev}_v(D^n\gamma) \in \mathfrak{h}_v
            \text{ for all }v \in V\text{ and }n \geq 0\}.
\end{equation}
It is stable under \(D\). In the one-loop example with the diagonal vertex
condition, it is precisely the smooth loop algebra
\(\mathcal{L}\mathfrak{g}\), as was discussed above.

The restriction of \(\omega\) to
\(\mathfrak{g}_{\infty}(\Gamma, \mathfrak{h})\) is again a two-cocycle.
A similar argument using compactly supported paths shows that it is
nontrivial whenever \(k_e \neq 0\) for at least one edge \(e\).

Let
\begin{equation}
    \widehat{\mathfrak{g}}_{\infty}(\Gamma, \mathfrak{h})
        = \mathfrak{g}_{\infty}(\Gamma, \mathfrak{h})
            \oplus \mathbb{C}K
\end{equation}
be the restriction of the central extension
\eqref{eq:central-extension} to this subalgebra.

\begin{proposition}
The derivation \(D\) of
\(\mathfrak{g}_{\infty}(\Gamma, \mathfrak{h})\) lifts to a derivation of
\(\widehat{\mathfrak{g}}_{\infty}(\Gamma, \mathfrak{h})\) by setting
\(D(K)=0\).
\end{proposition}

\begin{proof}
Since \(D\) is already a derivation on
\(\mathfrak{g}_{\infty}(\Gamma, \mathfrak{h})\), its extension with
\(D(K) = 0\) is a derivation precisely when the central terms in the
derivation identity cancel.
For \(\gamma, \eta \in
\mathfrak{g}_{\infty}(\Gamma, \mathfrak{h})\), integration by parts gives
\begin{align}
    \omega(D\gamma, \eta) + \omega(\gamma, D\eta)
        &= \sum_{e \in E} k_e
            \big(\langle \gamma_e(1),
                \partial_{\tau_e}\eta_e(1) \rangle_e
                - \langle \gamma_e(0),
                    \partial_{\tau_e}\eta_e(0) \rangle_e\big)
                    \notag\\
        &= \sum_{v \in V}
            \langle \gamma_v, (D\eta)_v \rangle_v = 0.
\end{align}
The last equality follows because both vertex values belong to
\(\mathfrak{h}_v\).
\end{proof}

The lifted derivation \(D\) defines an action of the one-dimensional
abelian Lie algebra \(\mathbb{C}d\). We use this action to construct a new
Lie algebra as the semidirect product
\begin{equation}\label{eq:derivation-extension}
    \widetilde{\mathfrak{g}}_{\infty}(\Gamma, \mathfrak{h})
        = \widehat{\mathfrak{g}}_{\infty}(\Gamma, \mathfrak{h})
            \rtimes \mathbb{C}d,
\end{equation}
where
\begin{equation}
    [d, \gamma] = D\gamma,
    \qquad
    [d, K] = 0.
\end{equation}
For nonzero edge weights, the central extension followed by this
derivation extension is the \emph{double-extension} construction, as in
affine Kac--Moody theory.

\subsection{An invariant bilinear form}
\label{subsec:invariant-form}

Fix the parametrizations \(\tau_e \in [0, 1]\) of all edges \(e \in E\).
Define a symmetric bilinear form on \(\mathfrak{g}(\Gamma, \mathfrak{h})\)
by
\begin{equation}
    B_0(\gamma, \eta)
        = \sum_{e \in E} k_e \int_0^1
            \langle \gamma_e(\tau_e),
                \eta_e(\tau_e) \rangle_e\,\mathrm{d}\tau_e.
\end{equation}
Its restriction to \(\mathfrak{g}_{\infty}(\Gamma, \mathfrak{h})\) will
also be denoted by \(B_0\).

\begin{proposition}\label{prop:integrated-form}
The form \(B_0\) is invariant on \(\mathfrak{g}(\Gamma, \mathfrak{h})\)
and on its subalgebra \(\mathfrak{g}_{\infty}(\Gamma, \mathfrak{h})\).
It is nondegenerate on both algebras if and only if \(k_e \neq 0\) for
every edge \(e \in E\).
\end{proposition}

\begin{proof}
Invariance on \(\mathfrak{g}(\Gamma, \mathfrak{h})\) follows by integrating
the invariance identity for the Killing form on each edge. Invariance
is preserved on restriction to a Lie subalgebra.

Suppose that \(k_e \neq 0\) for every edge, and let
\(\gamma \in \mathfrak{g}(\Gamma, \mathfrak{h})\) satisfy
\(B_0(\gamma, \eta) = 0\) for every
\(\eta \in \mathfrak{g}_{\infty}(\Gamma, \mathfrak{h})\).
Test against paths supported on one edge, of the form
\(\eta_e(\tau_e) = f(\tau_e)y\), where
\(f \in C_c^\infty((0, 1), \mathbb{C})\) and
\(y \in \mathfrak{g}_e\). All endpoint derivatives of these paths vanish,
so they belong to \(\mathfrak{g}_{\infty}(\Gamma, \mathfrak{h})\). Thus
\begin{equation}
    k_e \int_0^1 f(\tau_e)
        \langle \gamma_e(\tau_e), y \rangle_e\,\mathrm{d}\tau_e = 0
    \qquad\text{for all } f, y.
\end{equation}
Since \(k_e \neq 0\), arbitrary test functions and nondegeneracy of the
Killing form imply \(\gamma_e = 0\) on \((0, 1)\), hence on \([0, 1]\)
by continuity. Thus \(\gamma = 0\), which proves nondegeneracy on both
algebras.

Conversely, if \(k_e = 0\) for some edge, any nonzero smooth path
\(\gamma\) compactly supported in that edge belongs to
\(\mathfrak{g}_{\infty}(\Gamma, \mathfrak{h})\) and satisfies
\(B_0(\gamma, \eta) = 0\) for every
\(\eta \in \mathfrak{g}(\Gamma, \mathfrak{h})\). Hence \(B_0\) is
degenerate on both algebras.
\end{proof}

For \(\gamma, \eta \in
\mathfrak{g}_{\infty}(\Gamma, \mathfrak{h})\), the isotropy conditions
give
\begin{equation}
    B_0(D\gamma, \eta) + B_0(\gamma, D\eta) = 0,
    \qquad
    \omega(\gamma, \eta) = B_0(\gamma, D\eta).
\end{equation}
Extend the restricted form \(B_0\) to a symmetric bilinear form \(B\) on
the algebra in \eqref{eq:derivation-extension} by
\begin{equation}
    B(K, d) = -1,
    \qquad
    B(K, K) = B(d, d) = 0,
\end{equation}
and, for every
\(\gamma \in \mathfrak{g}_{\infty}(\Gamma, \mathfrak{h})\),
\begin{equation}
    B(K, \gamma) = B(d, \gamma) = 0.
\end{equation}

\begin{proposition}
The form \(B\) on
\(\widetilde{\mathfrak{g}}_{\infty}(\Gamma, \mathfrak{h})\) is invariant.
It is nondegenerate if and only if \(B_0\) is nondegenerate on
\(\mathfrak{g}_{\infty}(\Gamma, \mathfrak{h})\). These equivalent
conditions hold if and only if \(k_e \neq 0\) for every edge \(e \in E\).
\end{proposition}

\begin{proof}
By Proposition~\ref{prop:integrated-form}, \(B_0\) is invariant, which
proves invariance of \(B\) on triples of currents. The identities above
give
\begin{equation}
    B([d, \gamma], \eta) = B_0(D\gamma, \eta)
        = -\omega(\gamma, \eta) = B(d, [\gamma, \eta]^{\wedge}).
\end{equation}
The remaining cases follow from the defining pairings and the centrality
of \(K\). Thus \(B\) is invariant.

Let \(\gamma + aK + bd\) be such that
\begin{equation}
    B(\gamma + aK + bd, z) = 0
    \qquad\text{for all }
    z \in \widetilde{\mathfrak{g}}_{\infty}(\Gamma, \mathfrak{h}).
\end{equation}
Pairing with
\(K\) and \(d\) gives \(-b = 0\) and \(-a = 0\), respectively. The
remaining condition is
\begin{equation}
    B_0(\gamma, \eta) = 0
    \qquad\text{for all }
    \eta \in \mathfrak{g}_{\infty}(\Gamma, \mathfrak{h}).
\end{equation}
Conversely, any \(\gamma\) satisfying this condition also satisfies
\(B(\gamma, z) = 0\) for every
\(z \in \widetilde{\mathfrak{g}}_{\infty}(\Gamma, \mathfrak{h})\).
This proves the equivalence of nondegeneracy of \(B\) and \(B_0\).
The criterion in terms of the weights now follows from
Proposition~\ref{prop:integrated-form}.
\end{proof}

\section{Vector fields and the graph Virasoro algebra}
\label{sec:vector-fields}

\subsection{Reparametrizations and compatible vector fields}

The Lie algebras \(\mathfrak{g}(\Gamma)\) and
\(\mathfrak{g}(\Gamma, \mathfrak{h})\), together with the edge-integral
cocycle \(\omega\), are independent of orientation-preserving smooth
edge reparametrizations under the natural pullback identifications.
Hence the central extension
\(\widehat{\mathfrak{g}}(\Gamma, \mathfrak{h})\) is also independent of
these parametrizations. However, defining \(D\) and \(B_0\) requires
fixing the edge parametrizations.

The subalgebra \(\mathfrak{g}_{\infty}(\Gamma, \mathfrak{h})\) generally
depends on the parametrizations, but it is preserved by a special class
of reparametrizations \(\tau_e = \phi_e(\sigma_e)\). Namely, each
\(\phi_e \colon [0, 1] \to [0, 1]\) is an orientation-preserving smooth
diffeomorphism fixing the endpoints, and for every vertex \(v\) and
every \(r \geq 1\), the value \(\phi_e^{(r)}(i)\) is the same for all
\((e, i) \in \operatorname{Star}(v)\). For
\(\widetilde{\gamma}_e = \gamma_e \circ \phi_e\), the chain rule gives,
at \(i \in \{0, 1\}\),
\begin{align}
    \partial_{\sigma_e}\widetilde{\gamma}_e(i)
        &= \phi'_e(i)\partial_{\tau_e}\gamma_e(i), \\
    \partial_{\sigma_e}^{2}\widetilde{\gamma}_e(i)
        &= \big(\phi'_e(i)\big)^2
            \partial_{\tau_e}^{2}\gamma_e(i)
            + \phi''_e(i)\partial_{\tau_e}\gamma_e(i).
\end{align}
At each vertex, the scalar coefficients are the same on all incident
half-edges. Thus the new derivative evaluations are linear combinations
of the old ones in \(\mathfrak{h}_v\). The higher-order chain rule has
the same property, so all derivative vertex conditions are preserved.
The inverse reparametrizations have the same compatibility property.
This condition is sufficient; particular choices of \(\mathfrak{h}_v\)
may allow more reparametrizations.

The compatible reparametrizations above arise from flows of possibly
time-dependent real edgewise vector fields that vanish at the endpoints
and whose derivatives agree on all half-edges incident to each vertex.

The \emph{Lie algebra of graph vector fields}
\(\operatorname{Vec}(\Gamma)\) consists of the edgewise vector fields
\(X_f = (f_e\partial_{\tau_e})_{e \in E}\), where
\(f_e \in C^\infty([0, 1], \mathbb{C})\), satisfying the vertex conditions
\begin{equation}
    \partial_{\tau_e}^{r}f_e(i)
        = \partial_{\tau_{e'}}^{r}f_{e'}(j)
    \qquad
    \text{for }(e, i), (e', j) \in \operatorname{Star}(v),
    \quad v \in V,\quad r \geq 0.
\end{equation}
It is equipped with the standard Lie bracket
\begin{equation}
    [X_f, X_g]_e
        = \big(f_e\partial_{\tau_e}g_e
            - g_e\partial_{\tau_e}f_e\big)\partial_{\tau_e}.
\end{equation}

The Lie algebra \(\operatorname{Vec}(\Gamma)\) acts on
\(\mathfrak{g}_{\infty}(\Gamma, \mathfrak{h})\) by
\begin{equation}
    (D_{X_f}\gamma)_e = f_e\partial_{\tau_e}\gamma_e.
\end{equation}
Similarly to \(D\), the derivation \(D_{X_f}\) preserves the vertex
conditions in \(\mathfrak{g}_{\infty}(\Gamma, \mathfrak{h})\).
Writing \(f_v\) for the common value of
\(f_e(i)\) at \(v\), integration by parts gives
\begin{equation}
    \omega(D_{X_f}\gamma, \eta) + \omega(\gamma, D_{X_f}\eta)
        = \sum_{v \in V} f_v
            \langle \gamma_v, (D\eta)_v \rangle_v = 0.
\end{equation}
\begin{samepage}
Thus this action lifts to
\(\widehat{\mathfrak{g}}_{\infty}(\Gamma, \mathfrak{h})\) by setting
\(D_{X_f} K = 0\). Since
\([D_{X_f}, D_{X_g}] = D_{[X_f, X_g]}\), we can form the
semidirect-product extension
\begin{equation}
    \widehat{\mathfrak{g}}_{\infty}(\Gamma, \mathfrak{h})
        \rtimes \operatorname{Vec}(\Gamma),
\end{equation}
with mixed brackets \([X_f, \gamma] = D_{X_f}\gamma\) and
\([X_f, K] = 0\).
\end{samepage}

\subsection{Graph Virasoro algebra}

For the one-loop graph with the diagonal vertex condition, the algebra
\(\mathfrak{g}_{\infty}(\Gamma, \mathfrak{h})\) is the smooth loop
algebra. Infinitesimal circle reparametrizations act on it by vector
fields. The Fourier-polynomial vector fields form the Witt algebra.
Its standard central extension, defined by the Gelfand--Fuchs cocycle
\cite{GelfandFuchs}, is the Virasoro algebra.

\begin{proposition}\label{prop:graph-gelfand-fuchs}
Assume that every vertex \(v\) of \(\Gamma\) has equally many incoming
and outgoing half-edges, that is,
\(d_{\mathrm{in}}(v) = d_{\mathrm{out}}(v)\). Then the bilinear form
\begin{equation}
    \omega_{\mathrm{GF}}(X_f, X_g) = \sum_{e \in E} \int_0^1
        \big(g_e\partial_{\tau_e}^{3}f_e
            - f_e\partial_{\tau_e}^{3}g_e\big)\,\mathrm{d}\tau_e,
\end{equation}
for \(X_f, X_g \in \operatorname{Vec}(\Gamma)\), is a two-cocycle on
\(\operatorname{Vec}(\Gamma)\).
\end{proposition}

\begin{proof}
The form is skew-symmetric by definition.
For \(X_f, X_g, X_h \in \operatorname{Vec}(\Gamma)\), we compute the
cyclic sum explicitly on each edge. Differentiating the bracket
coefficient three times gives
\begin{equation}
    \partial_{\tau_e}^{3}\big(
        f_e\partial_{\tau_e}g_e - g_e\partial_{\tau_e}f_e\big)
        = f_e\partial_{\tau_e}^{4}g_e
            - g_e\partial_{\tau_e}^{4}f_e
        + 2\big(
            \partial_{\tau_e}f_e\,\partial_{\tau_e}^{3}g_e
            - \partial_{\tau_e}g_e\,\partial_{\tau_e}^{3}f_e\big).
\end{equation}
Here \(\sum_{\mathrm{cyc}}\) denotes summation over
\((f, g, h)\), \((g, h, f)\), and \((h, f, g)\). The fourth-derivative
terms cancel pairwise, and cyclic relabelling of the remaining terms
gives
\begin{align}
    \sum_{\mathrm{cyc}} h_e\big(
        f_e\partial_{\tau_e}^{4}g_e
        - g_e\partial_{\tau_e}^{4}f_e\big) &= 0, \\
    \sum_{\mathrm{cyc}} h_e\big(
        \partial_{\tau_e}f_e\,\partial_{\tau_e}^{3}g_e
        - \partial_{\tau_e}g_e\,\partial_{\tau_e}^{3}f_e\big)
        &= \sum_{\mathrm{cyc}}\big(
            f_e\partial_{\tau_e}g_e - g_e\partial_{\tau_e}f_e\big)
            \partial_{\tau_e}^{3}h_e.
\end{align}
Consequently, the cyclic sum of the edge integrands in
\(\sum_{\mathrm{cyc}} \omega_{\mathrm{GF}}([X_f, X_g], X_h)\) is
\begin{align}
    &\sum_{\mathrm{cyc}}\Big(
        h_e\partial_{\tau_e}^{3}\big(
            f_e\partial_{\tau_e}g_e - g_e\partial_{\tau_e}f_e\big)
        - \big(f_e\partial_{\tau_e}g_e
            - g_e\partial_{\tau_e}f_e\big)
            \partial_{\tau_e}^{3}h_e\Big) \\
    &\qquad = 2\sum_{\mathrm{cyc}}\big(
        f_e\partial_{\tau_e}g_e - g_e\partial_{\tau_e}f_e\big)
        \partial_{\tau_e}^{3}h_e
        - \sum_{\mathrm{cyc}}\big(
            f_e\partial_{\tau_e}g_e - g_e\partial_{\tau_e}f_e\big)
            \partial_{\tau_e}^{3}h_e \\
    &\qquad = \sum_{\mathrm{cyc}}\big(
        f_e\partial_{\tau_e}g_e - g_e\partial_{\tau_e}f_e\big)
        \partial_{\tau_e}^{3}h_e \\
    &\qquad = \det\begin{pmatrix}
        f_e & g_e & h_e \\
        \partial_{\tau_e}f_e & \partial_{\tau_e}g_e
            & \partial_{\tau_e}h_e \\
        \partial_{\tau_e}^{3}f_e & \partial_{\tau_e}^{3}g_e
            & \partial_{\tau_e}^{3}h_e
    \end{pmatrix}
    = \partial_{\tau_e}\det\begin{pmatrix}
        f_e & g_e & h_e \\
        \partial_{\tau_e}f_e & \partial_{\tau_e}g_e
            & \partial_{\tau_e}h_e \\
        \partial_{\tau_e}^{2}f_e & \partial_{\tau_e}^{2}g_e
            & \partial_{\tau_e}^{2}h_e
    \end{pmatrix}.
\end{align}
The determinant equality follows by expansion along the last row.
For the last equality, differentiating either of the first two rows in
the final determinant produces a determinant with two equal rows, so
only the derivative of the last row contributes.

The vertex conditions give a common endpoint value \(W_v\) of the final
determinant on all half-edges incident to \(v\). Integrating on each
edge and grouping the endpoint contributions by vertex gives
\begin{align}
    \sum_{\mathrm{cyc}} \omega_{\mathrm{GF}}([X_f, X_g], X_h)
        &= \sum_{e \in E} \int_0^1
            \partial_{\tau_e}\det\begin{pmatrix}
                f_e & g_e & h_e \\
                \partial_{\tau_e}f_e & \partial_{\tau_e}g_e
                    & \partial_{\tau_e}h_e \\
                \partial_{\tau_e}^{2}f_e & \partial_{\tau_e}^{2}g_e
                    & \partial_{\tau_e}^{2}h_e
            \end{pmatrix}\,\mathrm{d}\tau_e \\
        &= \sum_{v \in V}
            \big(d_{\mathrm{in}}(v) - d_{\mathrm{out}}(v)\big)W_v
            = 0.
\end{align}
The last equality uses the hypothesis at every vertex and proves the
cocycle identity.
\end{proof}

Under the assumption of Proposition~\ref{prop:graph-gelfand-fuchs},
define the \emph{graph Virasoro algebra} as the central extension
\begin{equation}
    \operatorname{Vir}(\Gamma)
        = \operatorname{Vec}(\Gamma) \oplus \mathbb{C}c
\end{equation}
with Lie brackets
\begin{equation}
    [X_f, X_g]_{\operatorname{Vir}}
        = [X_f, X_g] + \omega_{\mathrm{GF}}(X_f, X_g)c,
    \qquad [X_f, c]_{\operatorname{Vir}} = 0.
\end{equation}
Here the bracket on the right is the bracket in
\(\operatorname{Vec}(\Gamma)\), and \(c\) is a new central generator.

We can similarly form the semidirect-product extension
\begin{equation}
    \widehat{\mathfrak{g}}_{\infty}(\Gamma, \mathfrak{h})
        \rtimes \operatorname{Vir}(\Gamma)
        = \widetilde{\mathfrak{g}}^{\operatorname{Vir}}_{\infty}
            (\Gamma, \mathfrak{h}),
\end{equation}
with brackets
\begin{equation}
    [X_f, \gamma] = D_{X_f}\gamma,
    \qquad [X_f, K] = 0,
    \qquad [c, \gamma] = 0,
    \qquad [c, K] = 0.
\end{equation}


\begin{thebibliography}{99}

\bibitem{Borcherds}
R.~Borcherds, \emph{Monstrous moonshine and monstrous Lie superalgebras},
Invent.~Math. \textbf{109} (1992), 405--444.

\bibitem{DMS}
P.~Di~Francesco, P.~Mathieu, and D.~S\'en\'echal,
\emph{Conformal Field Theory}, Graduate Texts in Contemporary Physics,
Springer-Verlag, New York, 1997.

\bibitem{EF}
P.~I.~Etingof and I.~B.~Frenkel,
\emph{Central extensions of current groups in two dimensions},
Comm.~Math.~Phys. \textbf{165} (1994), no.~3, 429--444.

\bibitem{FF}
A.~Feingold and I.~Frenkel, \emph{A hyperbolic Kac--Moody algebra and the
theory of Siegel modular forms of genus 2},
Math.~Ann. \textbf{263} (1983), 87--144.

\bibitem{KM-int}
A.~J.~Feingold, A.~Kleinschmidt, and H.~Nicolai,
\emph{Hyperbolic Weyl groups and the four normed division algebras},
J.~Algebra \textbf{322} (2009), no.~4, 1295--1339; corrigendum,
J.~Algebra \textbf{489} (2017), 586--587.

\bibitem{GelfandFuchs}
I.~M.~Gel'fand and D.~B.~Fuchs,
\emph{The cohomologies of the Lie algebra of the vector fields in a
circle}, Funct.~Anal.~Appl. \textbf{2} (1968), no.~4, 342--343.

\bibitem{KacBook}
V.~G.~Kac, \emph{Infinite-Dimensional Lie Algebras}, 3rd ed.,
Cambridge University Press, Cambridge, 1990.

\bibitem{KZ}
V.~Knizhnik and A.~Zamolodchikov,
\emph{Current algebra and Wess--Zumino model in two dimensions},
Nucl.~Phys.~B \textbf{247} (1984), 83--103.

\bibitem{MO}
D.~Maulik and A.~Okounkov,
\emph{Quantum groups and quantum cohomology},
Ast\'erisque \textbf{408}, Soci\'et\'e Math\'ematique de France, Paris,
2019, ix+209 pp.

\end{thebibliography}
\end{document}